\documentclass[12pt,a4paper]{article}
\def\Dj{\hbox{D\kern-.73em\raise.30ex\hbox{-}
\raise-.30ex\hbox{}}}
\def\dj{\hbox{d\kern-.33em\raise.80ex\hbox{-}
\raise-.80ex\hbox{\kern-.40em}}}

\usepackage{amsmath}
\usepackage{amsthm}
\usepackage{amsfonts}
\usepackage{amssymb}
\usepackage{graphicx}
\allowdisplaybreaks[4]

\newtheorem{thm}{Theorem}
\newtheorem{cor}{Corollary}

\newtheorem{de}{Definition}

\newenvironment{Proof}
{\medskip\noindent\textbf{Proof.}}{\hfill$\square$ \par}

\begin{document}

\baselineskip=0.40in

\vspace*{30mm}

\begin{center}
{\bf \LARGE LAPLACIAN ESTRADA INDEX OF TREES}
\vspace{10mm}

{\large \bf Aleksandar Ili\' c$^a$,  Bo Zhou$^b$\footnote[3]{Corresponding author.}}

\vspace{6mm}

\baselineskip=0.20in

{\it $^a$Faculty of Sciences and Mathematics, University of
Ni\v s,\\
Vi\v segradska 33, 18000 Ni\v s, Serbia\/} \\
{\rm e-mail:} { \tt aleksandari@gmail.com}\\[2mm]

{\it $^b$Department of Mathematics, South China Normal University, \\
Guangzhou 510631, P. R. China\/} \\
{\rm e-mail:} {\tt zhoubo@scnu.edu.cn}

\vspace{6mm}

(Received August 25, 2009)
\end{center}

\vspace{6mm}

\baselineskip=0.25in

\noindent {\bf Abstract}

\vspace{4mm}

Let $G$ be a simple graph with $n$ vertices and let $\mu_1 \geqslant
\mu_2 \geqslant \ldots \geqslant \mu_{n - 1} \geqslant \mu_n = 0$ be
the eigenvalues of its Laplacian matrix. The Laplacian Estrada index
of a graph $G$ is defined as $LEE (G) = \sum\limits_{i = 1}^n e^{\mu_i}$.
Using the recent connection between Estrada index of a line graph
and Laplacian Estrada index, we prove that the path $P_n$ has
minimal, while the star $S_n$ has maximal $LEE$ among trees on $n$
vertices. In addition, we find the unique tree with the second
maximal Laplacian Estrada index.

\vspace{4mm}

\baselineskip=0.27in

\begin{center}
\section*{\normalsize 1. INTRODUCTION}
\end{center}

Let $G$ be a simple graph with $n$ vertices. The spectrum of $G$ consists of the eigenvalues
$\lambda_1 \geqslant \lambda_2 \geqslant \ldots \geqslant \lambda_n$ of its
adjacency matrix \cite{CvDS80}. The Estrada index of $G$ is defined as
\begin{equation}
\label{eq-estrada}
EE (G) = \sum_{i = 1}^n e^{\lambda_i}.
\end{equation}
This graph--spectrum--based graph invariant
was put forward by Estrada in \cite{Es00,Es02}, where it was shown that
$EE(G)$ can be used as a measure of the degree of folding of long
chain polymeric molecules. Further, it was shown in~\cite{Es07} that
the Estrada index provides a measure of the centrality of complex
networks, while a connection between the Estrada index and the
concept of extended atomic branching was pointed out
in~\cite{EsRR06}. Some mathematical properties of the Estrada index
were studied in~\cite{GuFMG07}--\cite{IlSt09}.

For a graph  $G$  with $n$ vertices, let $\mu_1 \geqslant \mu_2 \geqslant \ldots \geqslant \mu_{n - 1}
\geqslant \mu_n = 0$ be the eigenvalues of its Laplacian matrix \cite{Me94}.
In full analogy with Equation (\ref{eq-estrada}), the Laplacian
Estrada index of  $G$ is defined as~\cite{FTAG09}
$$
LEE (G) = \sum_{i = 1}^n e^{\mu_i}.
$$
Bounds for the Laplacian Estrada index may be found in
\cite{FTAG09}--\cite{ZhGu09}.
%\cite{FTAG09,Zh09,ZhGu09}.

\vspace{0.2cm}
Let $\mathcal {L} (G)$ be the line graph of $G$. In \cite{ZhGu09}, the authors proved the following
relation between Laplacian Estrada index of $G$ and Estrada index of a line graph of $G$.

\begin{thm} \textnormal{\cite{ZhGu09}}
\label{thm-zhou}
Let $G$ be a graph with $n$ vertices and $m$ edges. If $G$ is bipartite, then
$$
LEE (G) = n - m + e^2 \cdot EE (\mathcal{L} (G))\,.
$$
\end{thm}

Our goal here is to add some further evidence to support
the use of $LEE$ as a measure of branching in alkanes.
While the measure of branching cannot be formally defined,
there are several properties that any proposed measure has to satisfy \cite{MIPG,FiGuHRVV02}.
Basically, a topological index $TI$ acceptable as a measure of branching must satisfy the inequalities
$$
TI (P_n) < TI (T) < TI (S_n) \qquad \mbox{or} \qquad TI (P_n) > TI (T) > TI (S_n),
$$
for $n = 5, 6, \ldots$, where $P_n$ is the path,  $S_n$ is the star on $n$ vertices, and $T$ is any $n$-vertex tree, different from $P_n$ and $S_n$.
For example, the first relation is obeyed by the largest graph eigenvalue \cite{LoPe73} and Estrada index \cite{De09},
while the second relation is obeyed by the Wiener index \cite{EnJS76}, Hosoya index and graph energy \cite{Gu77}.
%\vspace{0.2cm}

We show that among the $n$-vertex trees, the path $P_n$ has minimal and
the star $S_n$ maximal Laplacian Estrada index,
$$
LEE(P_n) < LEE(T) < LEE(S_n),
$$
where $T$ is any $n$-vertex tree, different from $P_n$ and $S_n$. We
also find the unique tree with the second maximal Laplacian Estrada
index.

\begin{center}
\section*{\normalsize 2. LAPLACIAN ESTRADA INDEX OF TREES}
\end{center}

In our proofs, we will use a connection between Estrada index and the spectral moments of a graph.
For $k\geqslant 0$, we denote by~$M_k$ the $k$th spectral moment of~$G$,
$$
M_k=M_k (G) = \sum_{i = 1}^n \lambda_i^k.
$$
A walk of length~$k$ in~$G$ is any sequence of vertices and edges
of~$G$,
$$
w_{0},e_{1},w_{1},e_{2},\dots,w_{k - 1}, e_{k},w_{k},
$$
such that $e_{i}$ is the edge joining $w_{i-1}$ and~$w_{i}$ for
every $i=1,2, \ldots,k$. The walk is closed if $w_{0}=w_{k}$. It is
well-known (see~\cite{CvDS80}) that $M_k (G)$ represents the number
of closed walks of length~$k$ in~$G$. Obviously, for every graph
$M_0=n$, $M_1 = 0$ and $M_2 = 2m$. From the Taylor expansion of $e^x$, we have that the Estrada index
and the spectral moments of $G$ are related by
\begin{equation}
\label{eq-taylor} EE (G) = \sum_{k = 0}^{\infty} \frac{M_k}{k!}\,.
\end{equation}
Thus, if for two graphs $G$ and~$H$ we have $M_{k}(G)\geqslant
M_{k}(H)$ for all $k\geqslant 0$, then $EE(G)\geqslant EE(H)$.
Moreover, if the strict inequality $M_{k}(G)>M_{k}(H)$ holds for at
least one value of~$k$, then $EE(G)>EE(H)$.
%\vspace{0.2cm}

Among the $n$-vertex connected graphs, the path $P_n$ has minimal and the complete graph $K_n$ maximal Estrada index \cite{GuER07,De09},
\begin{equation}
\label{eq-maxmin}
EE(P_n) < EE(G) < EE(K_n),
\end{equation}
where $G$ is any $n$-vertex connected graph, different from $P_n$ and $K_n$.

\begin{thm}
Among the $n$-vertex trees, the path $P_n$ has minimal and
the star $S_n$ maximal Laplacian Estrada index,
$$
LEE(P_n) < LEE(T) < LEE(S_n),
$$
where $T$ is any $n$-vertex tree, different from $P_n$ and $S_n$.
\end{thm}

\begin{Proof}
The line graph of a tree $T$ is a connected graph with $n - 1$
vertices. The line graph of a path $P_n$ is also a path  $P_{n-1}$,
while the line graph of a star $S_n$ is a complete graph $K_{n-1}$.
Using the relation (\ref{eq-maxmin}) it follows that
$$
EE (\mathcal{L} (P_n)) \leqslant EE (\mathcal{L} (T)) \leqslant EE
(\mathcal{ L} (S_n)),
$$
and from Theorem \ref{thm-zhou} we get $LEE(P_n) \leqslant LEE(T)
\leqslant LEE(S_n)$ with left equality if and only if $T \cong P_n$
and right equality if and only if $T \cong S_n$.
\end{Proof}

\begin{center}
\section*{\normalsize 3. SECOND MAXIMAL LAPLACIAN ESTRADA OF TREES}
\end{center}

\begin{de}
Let $v$ be a vertex of degree $p + 1$ in a graph $G$, which is not a
star, such that $v v_1, v v_2, \ldots, v v_p$ are pendent edges
incident with $v$ and $u$ is the neighbor of $v$ distinct from $v_1,
v_2, \ldots,$ $v_p$. We form a graph $G' = \sigma (G, v)$ by removing
edges $v v_1, v v_2, \dots, v v_p$ and adding new edges $u v_1, u
v_2, \dots,$ $u v_p$. We say that $G'$ is $\sigma$-transform of $G$.
\end{de}

\begin{figure}[ht]
  \center
  \includegraphics [width = 11cm]{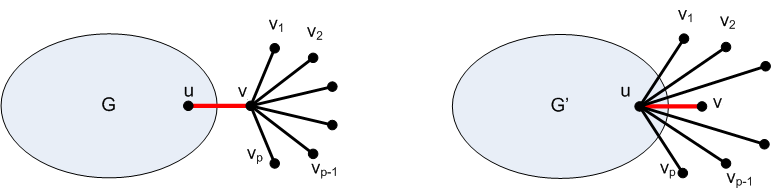}\vspace{0.2cm}
  \caption { $\sigma$-transformation applied to $G$ at vertex $v$. }
\end{figure}

%\vspace{0.2cm}

\begin{thm}
\label{sigma} Let $G' = \sigma (G, v)$ be a $\sigma$-transform of a
bipartite graph $G$. Then
\begin{equation}
\label{eq-sigma} LEE (G) < LEE (G')\,.
\end{equation}
\end{thm}

\begin{Proof}
The graphs $G$ and $G'$ are both bipartite and have the same number of vertices and edges.
Using Theorem \ref{thm-zhou}, it is enough to prove inequality
$$
EE (\mathcal{L} (G)) < EE (\mathcal{L} (G'))\,.
$$
Let $u_1, u_2, \dots, u_m$ be the neighbors of $u$ in $G$, different from $v$.
Consider the induced subgraph $H$ of $\mathcal{L} (G)$ formed using vertices
$v v_1, v v_2, \dots, v v_p$, $v u$, $u u_1, u u_2, \dots, u u_m$. It is easy to see that
these vertices are grouped in two cliques of sizes $p + 1$ and $m + 1$ with the vertex $uv$ in common.
Similarly, consider the induced subgraph $H'$ of $\mathcal{L} (G')$ formed using corresponding vertices
$u v_1, u v_2, \dots, u v_p$, $v u$, $u u_1, u u_2, \dots, u u_m$. Here, we have one clique of size $m + p + 1$.

Since $H$ is a proper subgraph of $H'$, it follows that for every $k \geqslant 0$, $M_{k}(H')\geqslant
M_{k}(H)$ and then $M_k({\mathcal L} (G'))\ge M_k({\mathcal L} (G))$, which is strict for some $k$. Finally, using the relation (\ref{eq-taylor}), we have
$LEE (G) < LEE (G')$.
\end{Proof}

%\begin{figure}[ht]
%  \center
%  \includegraphics [width = 11cm]{sigma_transformation.png}
%  \caption { $\sigma$-transformation applied to $G$ at vertex $v$. }
%\end{figure}
%
%\vspace{0.2cm}

Let $T$ be an arbitrary tree on $n$ vertices with root $v$. We can
find a vertex $u$ that is the parent of the leaf on the deepest
level and apply $\sigma$-transformation at $u$ to strictly increase the
Laplacian Estrada index.

\begin{cor}
Let $T$ be a tree on $n$ vertices. If $T
\not \cong S_n$, then $LEE (T) < LEE (S_n)$.
\end{cor}

Let $S_n (a, b)$ be the tree formed by adding an edge between the centers of the stars
$S_a$ and $S_b$, where $a + b = n$ and $2\leqslant a\leqslant \lfloor \frac{n}{2} \rfloor$. We call $S_n (a, b)$ the double star.
By direct calculation, the characteristic polynomial of the Laplacian matrix of the double star $S_n (a, b)$ is equal to
$$
P (x) = (-1)^n x (x-1)^{n-4} \left( x^3-(n+2)x^2+(n+2+ab)x-n \right ).
$$

We may assume that %$1 \leqslant a \leqslant\frac{n}{2}$ and
$n > 5$.
The Laplacian spectra of $S_n (a, b)$ consists of three real roots
of polynomial $f_{n, a} (x) = x^3-(n+2)x^2+(n+2+a (n - a))x-n$, $1$
with multiplicity $n - 4$, and $0$ with multiplicity one. In order
to establish the ordering of double stars with $n$ vertices by $LEE$
values it is enough to consider the following function
$
g_{n, a} (x_1, x_2, x_3) = e^{x_1} + e^{x_2} + e^{x_3}$,
where $x_1 \geqslant x_2 \geqslant x_3 > 0$ are the roots of $f_{n,a} (x)$.
%From Vieta's Formulas, we have $x_1 + x_2 + x_3 = n + 2$.

%For the largest root $x_1$, we have
We locate $x_1$, $x_2$ and $x_3$. First we have
$$
f_{n, a} (n - a + 1) = 1 - a < 0
$$
and
$$
f_{n, a} \left( n - a + \frac{3}{2} \right) = \frac{15}{8} + a^2 + n + \frac{n^2}{2} - \frac{11a}{4} - \frac{3na}{2}\,.
$$
The last function (considered as a quadratic function of $a$) is decreasing for $a < \frac{11}{8} + \frac{3n}{4}$,
and then for $a \leqslant  \frac{n}{2} - 1$, we have
$$
f_{n, a} \left( n - a + \frac{3}{2} \right) \geqslant f_{n, a} \left( n - \frac{n}{2} + 1 + \frac{3}{2} \right) = \frac{45 + n}{8} > 0\,.
$$
%%with the special case $a = \lfloor \frac{n}{2} \rfloor$.
%For the middle root $x_2$, we have
Next we have
$$
f_{n, a} (a) = (a - 1) (n - 2 a) \geqslant 0 \ \mbox{ and } \ f_{n, a} ( a + 1 ) = 1 + a - n < 0\,.
$$
%$$
%f_{n, a} (a) = (a - 1) (n - 2 a) \geqslant 0
%$$
%and
%$$
%f_{n, a} ( a + 1 ) = 1 + a - n < 0.
%$$
%For the smallest root $x_3$, we have
Finally we have
$$
f_{n, a} (0) = - n < 0 \ \mbox{ and } \  f_{n, a} (1) = (a -1) (n - 1 - a) > 0\,.
$$
%
%
%$$
%f_{n, a} (0) = - n < 0
%$$
%and
%$$
%%f_{n, a} \left( \frac{1}{2} \right) = \frac{4an - 4a^2 - 6n + 5}{8} > 0.
%f_{n, a} (1) = (a -1) (n - 1 - a) > 0.
%$$
Thus %$f_{n, a}(x)$ has roots in the following intervals
$x_3\in [0, 1]$, $x_2\in [a, a + 1]$ and $x_1\in \left[n - a + 1, n - a + \frac{3}{2}\right]$ for $2\leqslant a \leqslant  \frac{n}{2} - 1$.
The function
$$
h (a) = e^0+ e^a + e^{n-a+1} -e^1 - e^{a+2} - e^{n-a+1/2}
$$
is decreasing for $a > 0$ (since $h' (a) < 0$), and then for $a \leqslant \frac{n}{2} - 1$ we have
$$
h(a) \geqslant h \left ( \frac{n}{2} - 1 \right) = e^{n/2} \left ( e^{-1} - e + e^2 - e^{3/2} \right) +1-e > \frac{e^{n/2}}{2} +1-e > 0\,.
$$
Thus, for $2 \leqslant a < \lfloor \frac{n}{2} \rfloor -1$, we have
$$
e^0 + e^a + e^{n-a+1} > e^1 + e^{a+2} + e^{n-a+1/2},
$$
and then
$LEE (S_n (a, b)) > LEE (S_n (a + 1, b - 1))$.
The special case $a = \lfloor \frac{n}{2} \rfloor$ can be handled easily,
\begin{eqnarray*}
&&LEE (S_n (2, n - 2)) - LEE \left(S_n \left(\left\lfloor \frac{n}{2}
\right\rfloor, \left\lceil \frac{n}{2} \right\rceil\right)\right)\\[1.5mm]
&>& e^{n-1}+e^2 - e^{\lceil n/2 \rceil} - e^{\lceil n/2 \rceil + 2} - e %\\
%&>&
>e^{\lceil n/2 \rceil} \cdot \left ( e^{\lfloor n/2 \rfloor -1} - 1 - e^2 \right ) - e > 0 \mbox{ for } n>7
\end{eqnarray*}
and  by direct calculation, we also have
$LEE (S_n (2, n - 2)) - LEE \left(S_n \left(\left\lfloor \frac{n}{2}
\right\rfloor, \left\lceil \frac{n}{2} \right\rceil\right)\right)>0$ for $n=6, 7$.
By Theorem \ref{sigma}, the second maximal $LEE$ for $n$-vertex trees
is a double star $S_n(a,b)$, and then from discussions above, we have

\begin{cor}
The unique tree on $n \geqslant 5$ vertices with the second maximal Laplacian Estrada index is a double star
$S_n (2, n - 2)$.
\end{cor}

Note that as above,  $LEE (S_n (3, n - 3)) - LEE \left(S_n \left(\left\lfloor \frac{n}{2}
\right\rfloor, \left\lceil \frac{n}{2} \right\rceil\right)\right)>0$ for $n\geqslant 8$ (the cases
for $n=8, 9$ need direct calculation).
%By Theorem \ref{sigma} and discussions above,  the third maximal Laplacian Estrada index for trees
%on $n\geqslant 6$ vertices is uniquely achieved by $S_n(3, n-3)$ or the tree $C_n(n-5)$ formed by
%attaching $n-5$ vertices to the center of a path $P_5$.
Tested by computer on trees with at most $22$ vertices, $S_n (3, n - 3)$ is the unique tree with
the third and $C_n (n-5)$ (formed by attaching $n-5$ vertices to the center of a path $P_5$) is the
unique tree with the fourth maximal Laplacian Estrada index for $n\geqslant 6$.

\baselineskip=0.25in

\bigskip
\noindent
{\bf Acknowledgement. } This work is supported by Research Grant 144007 of Serbian
Ministry of Science and Technological Development, and the Guangdong Provincial
Natural Science Foundation of China (Grant No. 8151063101000026).


\begin{thebibliography}{99}

\bibitem{CvDS80}
    D. Cvetkovi\' c, M. Doob, H. Sachs,
    \textit{Spectra of Graphs -- Theory and Application},
    Academic Press, New York, 1980.

\bibitem{Es00}
    E. Estrada,
    Characterization of 3D molecular structure,
    \textit{Chem. Phys. Lett.} \textbf{319} (2000) 713--718.

\bibitem{Es02}
    E. Estrada,
    Characterization of the folding degree of proteins,
    \textit{Bioinformatics} \textbf{18} (2002) 697--704.

\bibitem{Es07}
    E. Estrada,
    Topological structural classes of complex networks,
    \textit{Phys. Rev. E} \textbf{75} (2007) 016103-1--016103-12.

\bibitem{EsRR06}
    E. Estrada, J. A. Rodr\'{\i}guez-Val\'{a}zquez, M. Randi\' c,
    Atomic branching in molecules,
    \textit{Int. J. Quantum Chem.} \textbf{106} (2006) 823--832.

\bibitem{GuFMG07}
    I. Gutman, B. Furtula, V. Markovi\' c, B. Gli\v si\' c,
    Alkanes with greatest Estrada index,
    \textit{Z. Naturforsch.} \textbf{62a} (2007) 495--498.

\bibitem{GuGr07}
    I. Gutman, A. Graovac,
    Estrada index of cycles and paths,
    \textit{Chem. Phys. Lett.} \textbf{436} (2007) 294--296.

\bibitem{GuER07}
    I. Gutman, E. Estrada, J. A. Rodr\'{\i}guez-Val\'{a}zquez,
    On a graph--spectrum--based structure descriptor,
    \textit{Croat. Chem. Acta} \textbf{80} (2007) 151--154.

\bibitem{GuRGP07}
    I. Gutman, S. Radenkovi\' c, A. Graovac, D. Plav\v si\' c,
    Monte Carlo approach to Estrada index,
    \textit{Chem. Phys. Lett.} \textbf{447} (2007) 233--236.

\bibitem{PeGR07}
    J. A. de la Pe\~ na, I. Gutman, J. Rada,
    Estimating the Estrada index,
    \textit{Linear Algebra Appl.} \textbf{427} (2007) 70--76.

\bibitem{Gut08} I. Gutman, Lower bounds for Estrada index,
     \textit{Publ. Inst. Math. $($Beograd$)$ $($N.S.$)$} \textbf{83} (2008) 1--7.

\bibitem{Zh08}
    B. Zhou,
    On Estrada index,
    \textit{MATCH Commun. Math. Comput. Chem.} \textbf{60} (2008) 485--492.

\bibitem{RoJiMe09}
    M. Robbiano, R. Jim\' enez, L. Medina,
    The energy and an approximation to Estrada index of some trees,
    \textit{MATCH Commun. Math. Comput. Chem.} \textbf{61} (2009) 369--382.

\bibitem{De09}
    H. Deng,
    A proof of a conjecture on the Estrada index,
    \textit{MATCH Commun. Math. Comput. Chem.} \textbf{62} (2009) 599--606.

\bibitem{De09a}
    H. Deng,
    A note on the Estrada index of trees,
    \textit{MATCH Commun. Math. Comput. Chem.} \textbf{62} (2009) 607--610.

\bibitem{IlSt09}
    A. Ili\' c, D. Stevanovi\' c,
    The Estrada index of chemical trees,
    \textit{J. Math. Chem.} (2009), doi:10.1007/s10910-009-9570-0.

\bibitem{Me94} R. Merris, Laplacian matrices of graphs: a survey,
       \textit{Linear Algebra Appl.} \textbf{197--198} (1994) 143--176.

\bibitem{FTAG09}
    G. H. Fath-Tabar, A. R. Ashrafi, I. Gutman,
    Note on Estrada and $L$-Estrada indices of graphs,
    \textit{Bull. Cl. Sci. Math. Nat. Sci. Math.}, in press.

\bibitem{Zh09}
    B. Zhou,
    On sum of powers of Laplacian eigenvalues and Laplacian Estrada index of graphs,
    \textit{MATCH Commun. Math. Comput. Chem.} \textbf{62} (2009) 611--619.

\bibitem{ZhGu09}
    B. Zhou, I. Gutman,
    More on the Laplacian Estrada index,
    \textit{Appl. Anal. Discrete Math.} \textbf{3} (2009) 371--378.

\bibitem{MIPG}
    Z.~Markovi\'c, V.~Ivanov-Petrovi\'c, I.~Gutman,
    Extremely branched alkanes,
    \textit{J. Mol. Struct. $($Theochem$)$} \textbf{629} (2003) 303--306.

\bibitem{FiGuHRVV02}
    M. Fischermann, I. Gutman, A. Hoffmann, D. Rautenbach, D. Vidovi\' c, L. Volkmann,
    Extremal chemical trees,
    \textit{Z. Naturforsch.} \textbf{57a} (2002) 49--52.

\bibitem{LoPe73}
    L. Lov\'{a}sz, J. Pelik\'{a}n,
    On the eigenvalues of trees,
    \textit{Period. Math. Hung.} \textbf{3} (1973) 175--182.

\bibitem{EnJS76}
    R. C. Entringer, D. E. Jackson, D. A. Snyder,
    Distance in graphs,
    \textit{Czech. Math. J.} \textbf{26} (1976) 283--296.

\bibitem{Gu77}
    I. Gutman,
    Acyclic systems with extremal H\" uckel $\pi$-electron energy,
    \textit{Theor. Chim. Acta}  \textbf{45} (1977) 79--87.




\end{thebibliography}
\end{document}